\documentclass[11pt]{article}
\usepackage[english]{babel}
 \usepackage{makeidx}
 \usepackage{amssymb,amsbsy,amsmath,amsfonts}

\usepackage{graphics,graphicx}
\usepackage{color}
\def\v3{\vskip0.3cm \noindent}

\newtheorem{rem}{\bf Remark}[section]

\newtheorem{thm}{\bf Theorem}

\newtheorem{lem}{\bf Lemma\/}[section]

\def\QQ{{\rlap {\raise 0.4ex \hbox{$\scriptscriptstyle  $}}
\hskip -0.2em Q}}

\def\CC{{ \rlap {\raise 0.4ex \hbox{$\scriptscriptstyle  $}}
\hskip -0.2em C}}

\def\v3{\vskip0.3cm \noindent}
\date{today}
\begin{document}
\title{\Large \sf Estimate of the validity interval for the Antimaximum Principle and \\
application to a non-cooperative system} \author{\small \sf
J. FLECKINGER  \\
 {\small Institut de Math\'ematique - CEREMATH-UT1}  \\ {\small Universit\' e de Toulouse, } {\small 31042  Toulouse Cedex, France }\\ {\small
jfleckinger@gmail.com} \\ \\
{\small \sf J. HERNANDEZ}\\ {\small Departamento de Matem\'aticas  }\\ {\small Universidad Aut\'onoma,} {\small 28049
Madrid, Spain}\\ {\small  jesus.hernandez@uam.es}\\ \\
{\small \sf
F. de TH\' ELIN} \\ {\small Institut de Math\'ematique}\\{\small Universit\' e de Toulouse, } {\small 31062  Toulouse C\' edex, France \ }\\
{\small francoisdethelin@yahoo.fr}}

\date{\today}

\maketitle 
\vspace{-0.20cm}
\begin{abstract}
{\small We are concerned with the sign of the solutions of  non-cooperative systems  when the parameter varies near 
a principal  eigenvalue of the system. With this aim we give  precise estimates of 
 the validity interval for the Antimaximum Principle for an equation and 
an example. We apply these results to a non-cooperative system. Finally  a counterexample shows that our hypotheses are necessary.  The  Maximum Principle remains true only  for a restricted positive cone. }

 \end{abstract}
\newpage \noindent

\section{Introduction} 

 In this paper we use ideas concerning the Anti-Maximum Principle due to Cl\'ement and Peletier \cite{ClPe} and later to Arcoya Gamez \cite{ArGa}
to obtain in Section 2 precise estimates concerning the validity interval for the Anti-maximum Principle for 
one equation. An example shows that this estimate is sharp.
\par \noindent
The Maximum Principle and then the Anti-Maximum Principle 
for the case of a single equation have been extensively studied 
later for cooperative elliptic systems (see the references 
(\cite{Am},\cite{dFMi1986},\cite{FiMi1},\cite{FiMi2},\cite{FHeT},\cite{PW}). The results in \cite{FHeT}, 
are still valid for systems(with constant coefficients)  involving the $p$-Laplacian.
Some results for non-cooperative systems can be found  $e.g.$ in \cite{CaMi},\cite{Le}.
Very general results concerning the Maximum Principle for equations and cooperative
 systems for different classes (classical, weak, very weak) of solutions were given by Amann 
 in a long paper \cite{Am2005}, in particular the Maximum Principle was shown to be equivalent to the 
 positivity of the principal eigenvalue.
\par \noindent
Here in Section 3, we consider a non-cooperative $2 \times 2$ system with constant coefficients 
depending on a real parameter $\mu$ having two real principal eigenvalues $\mu_1^- < \mu_1^+$. We obtain 
some theorems of Anti-Maximum principle type concerning the behavior of different cones of couples of functions  having positivity (or negativity) properties. We give several results of this type for values of $\mu_1^- < \mu$ but 
close to $\mu_1^-$ by combining the usual Maximum Principle  and the results for the Anti-Maximum Principle 
in Section 2.
\par \noindent
 Finally a counterexample is given showing that the Maximum Principle does not hold in general
for non cooperative systems, but a (partial, under an additional assumption) Maximum  Principle for $\mu < \mu_1^-$ is also obtained. 

\section{Estimate of the validity interval for the anti-maximum principle} 
Let $\Omega$ be a smooth bounded domain in $I\!\!R^N$. We consider the following Dirichlet boundary value problem
\begin{equation} \label{1}
-\Delta z \, =\, \mu z + \,h \;\, {\rm in }\; \Omega \; ,  \;\; z=0 \; {\rm on } \; \partial \Omega,  \end{equation}
where $\mu$ is a real parameter. 
We associate to (\ref{1}) the eigenvalue problem 
\begin{equation} \label{2 }-\Delta \varphi \, =\, \lambda \varphi  \;\, {\rm in }\; \Omega \; , \;\; \varphi=0 \; {\rm on } \; \partial \Omega. \end{equation}
We denote by $\lambda_k$, $k\in I\!\!N^*$ the eigenvalues ($0<\lambda_1 < \lambda_2 \leq ...$) and by $\varphi_k$ a set of orthonormal associated eigenfunctions.
We choose $\varphi_1 >0$.
\par \noindent 
{\bf Hypothesis $(H_0)$:} $\;$ We write 
\begin{equation} \label{3}h=\alpha \varphi_1 + h^{\bot} \end{equation}
where $\int_{\Omega} h^{\perp} \varphi_1 =0$  and we assume $\alpha>0$ and $h\in L^q$, $q>N$ if $N\geq 2$ and $q=2$ if $N=1$. 
\begin{thm}: $\,$ 
We assume $(H_0)$ and $\lambda_1<\mu \leq  \Lambda < \lambda_2$. There exists a constant $K$ depending only on $\Omega$, $\Lambda$ and $q$ such that,
 for $\lambda_1<\mu< \lambda_1 + \delta(h)$ with
\begin{equation} \label{ E}\delta(h)=\frac{K \alpha}{\|h^{\bot}\|_{L^q}}, \end{equation}
the solution $z$ to (\ref{1}) satisfies the antimaximum principle,  that is
\begin{equation} \label{AMP }z<0 \;  {\rm in }\; \Omega; \;\; \partial z/ \partial \nu >0 \; {\rm on } \; \partial \Omega,  \end{equation}
where $\partial / \partial \nu$ denotes the outward normal derivative. 
\end{thm}
\begin{rem}
The antimaximum principle of Theorem 1, assuming $\alpha>0$, is in the line of the version given by Arcoya Gamez \cite{ArGa}.
\end{rem}

\begin{lem} \label{L21} We assume $ \lambda_1< \mu \leq \Lambda <\lambda_2$ and $h \in L^q, \, q>N\geq 2$.
We suppose that there exists a  constant $C_1$ depending only on $\Omega, q,$ and $\Lambda$  such that $z$ satisfying 
(\ref{1}) is such that
\begin{equation} \label{4} \|z\|_{L^2} \leq C_1 \|h\|_{L^2}.\end{equation}
Then there exist  constants $C_2$ and $C_3$, depending only on $\Omega, q$ and $\Lambda$ such that
\begin{equation} \label{5} \|z\|_{{\cal C}^1} \leq C_2 \|h\|_{L^q} \; {\rm and  } \; \| z \|_{L^q} \leq \, C_3 \|h\|_{L^q}. \end{equation}
\end{lem}

\begin{rem}
 Hypothesis \ref{4} cannot hold, unless $h$ is orthonal to $\varphi_1$. Indeed, letting $\mu$ go to $\lambda_1$, 
\ref{4} implies the existence of a solution to  \ref{1} with $\mu=\lambda_1$. Note that in the proof of Theorem 1,
Lemma \ref{L21} is used for $h$ orthogonal to $\varphi_1$.
\end{rem}

\subsection{ Proof of Lemma \ref{L21}}
All constants in this proof  depend only on $\Omega$, $\Lambda$ and   $q$.
\vskip0.15cm \noindent
{\bf Claim:} $\,$ $\| z \|_{L^q} \leq \, C_3 \|h\|_{L^q}$.
\par \noindent 
If the claim is verified then,  
by regularity results for the Laplace operator combined with Sobolev imbeddings 
\begin{equation} \label{6 }\| z\|_{{\cal C}^1} \, \leq \, C_4 \| z \|_{W^{2,q}} \, \leq \, 
C_5 ( \Lambda \| z \|_{L^q} +  \| h \|_{L^q} ).  \end{equation}
From the claim and   regularity results  we deduce (\ref{5}).
\vskip0.15cm \noindent
{\bf Proof of the claim: }
\par \noindent
{\bf - Step 1} $\;$ We consider the sequence $p_j= 2 + \frac{8j}{N}$ for $j\in I\!\!N$.
Observe that for any $j$, $W^{2, p_j}  \hookrightarrow L^{p_{j+1}}$ and that there exists a constant 
$H(j)$ such that 
\begin{equation}\label{7} \forall v \in W^{2,p_j}, \, \|v\|_{L^{p_{j+1}}} \; \leq \; H(j) \|v\|_{W^{2, p_j}}. \end{equation}
The relation (\ref{7}) is obvious if $2p_j\geq N$ and for $2p_j<N$ we have
$$\frac{Np_j}{N-2p_j} - p_{j+1} \; =  \; \frac{2p_j p_{j+1} - 8}{N-2p_j} \; > \; 0$$ and the result follows by classical 
Sobolev imbedding. 
\vskip0.15cm  \noindent
{\bf - Step 2} $\;$ We consider $z$ satisfying (\ref{1}). For $j=0$, we derive from (\ref{4}) and 
H\"older inequality that 
\begin{equation}\label{8}\|z\|_{L^2} \leq C_5 \|h\|_{L^q}.\end{equation}
By induction we assume that $z \in L^{p_{j}}$ with $p_j<q$ and that 
\begin{equation}\label{9} \|z\|_{L^{p_j}} \leq K(j)  \|h\|_{L^q}.\end{equation}
By H\"older inequality, 
$$ \|\mu z + h\|_{L^{p_j}} \leq \Lambda
 \|z\|_{L^{p_j}} + |\Omega|^{\frac{q-p_j}{q p_j}} \|h\|_{L^q} .$$ 
By regularity results for the Laplace operator:
$$\|z\|_{W^{2,p_j}} \;\leq\; 
C(j) ( \Lambda  \|z\|_{L^{p_j}} + |\Omega|^{\frac{q-p_j}{q p_j}} \|h\|_{L^q}) \, \leq \, C(j) 
(\Lambda K(j) +  |\Omega|^{\frac{q-p_j}{q p_j}} )  \|h\|_{L^q} .$$ 
Using (\ref{7})  the relation (\ref{9}) holds for $j+1$ and the induction is proved. 
\vskip0.15cm \noindent
{\bf - Step 3} $\;$ Let $J$ be such that $p_{J+1} \geq q > p_J$. After $J$ iterations we get
by (\ref{9})
$$\|z\|_{L^q} \, \leq \, C_6 \|z\|_{L^{p_{J+1}}} \, \leq \, C_6 K(J+1)\|z\|_{W^{2,p}}
\leq $$
$$ C_7 K(J+1) \|\mu z + h\|_{L^{p_{J}}} \leq C_8 (\Lambda \|h\|_{L^q} + \|h\|_{L^{p_{J}}} ) \leq 
C_9 \|h\|_{L^q}, $$
which is  the claim.
$\bullet$

\subsection{ Proof of Theorem 1}
 \noindent
{\bf - Step 1:} $\;$ We prove the following inequality: 
\begin{equation}\label{10}\| z^{\bot}\|_{{\cal C}^1} \, \leq \, C_2   \| h^{\bot}\|_{L^q}. \end{equation}
We derive from (\ref{3}) 
\begin{equation}\label{11}z\, =\, \frac{\alpha}{\lambda_1 - \mu} \varphi_1 \, + \, z^{\bot}, \end{equation}
with $z^{\bot}$ solution of 
\begin{equation}\label{12}  -\Delta z^{\bot} \, =\, \mu z^{\bot} + \,h^{\bot} \;\, {\rm in }\; \Omega \; ; \;\;
 z^{\bot}=0 \; {\rm on } \; \partial \Omega.  \end{equation}
By the  variational characterization of $\lambda_2$:
$$ \lambda_2 \int_{\Omega} |z^{\bot}|^2 \, \leq \,  \int_{\Omega} |\nabla z^{\bot}|^2 \, = 
\, \mu \int_{\Omega} |z^{\bot}|^2 \, +\, \int_{\Omega} z^{\bot} h^{\bot}.$$ 
Hence 
$$\| z^{\bot}\|_{L^2} \, \leq \, \frac{1}{\lambda_2-\Lambda} \| h^{\bot}\|_{L^2} .$$
By Lemma \ref{L21}, we derive (\ref{10}). 
\vskip0.15cm \noindent
{\bf - Step 2:} $\;$ {\it Close to the boundary:} 
\par \noindent
We show now that  on   the boundary
$\frac{\partial z }{  \partial \nu }(x) > 0.$ and near the boundary 
$z<0$.
\par \noindent
Since $\partial \varphi_1/ \partial \nu < 0 \; {\rm on } \; \partial \Omega$, we set 
\begin{equation}\label{13} A:= min_{\partial \Omega} |\partial \varphi_1/ \partial \nu | >0. \end{equation}
By a  continuity argument there exists $\varepsilon >0$ such that
\begin{equation}\label{14}  dist(x,  \partial \Omega) <\varepsilon \; \Rightarrow \; \partial \varphi_1/ \partial \nu  (x) \leq -A/2. \end{equation}
Hence by (\ref{10}) to (\ref{14}) , for any $x\in \Omega$ such that $dist(x,  \partial \Omega) <\varepsilon$, and if 
$$0< \mu - \lambda_1 < \frac{\alpha A} {4 C_2  \|h^{\bot}\|_{L^q}}, \; $$
we have 
$$\frac{\partial z }{  \partial \nu } (x) = \frac{\alpha}{\lambda_1 - \mu} \frac{\partial \varphi_1 }{  \partial \nu } (x) 
\, +\, \frac{\partial z^{\bot} }{  \partial \nu } (x)\, \geq \, 
\frac{\alpha}{\lambda_1 - \mu} \frac{\partial \varphi_1 }{  \partial \nu } (x) - C_2 \|h^{\bot}\|_{L^q} ,  $$
hence 
\begin{equation}\label{15} \frac{\partial z }{  \partial \nu } (x)  \geq  \; \frac{\alpha}{2(\lambda_1 - \mu)} \frac{\partial \varphi_1 }{  \partial \nu } (x) \, >0. 
\end{equation}
Therefore 
$\frac{\partial z }{  \partial \nu } (x) >0 \; {\rm on }\; \partial \Omega$. Moreover 
since $z=\varphi_1=0$ on $\partial \Omega$, we deduce from (\ref{15}) that, 
for $x\in \Omega$ with 
$ dist(x,  \partial \Omega) <\varepsilon' \leq \varepsilon/2   $ ($\varepsilon'$ small enough), 
$$z(x) \leq  \; \frac{\alpha}{2(\lambda_1 - \mu)} \varphi_1  (x) <0,  
$$
where $\varepsilon'$ does not depend on $\mu$. 
\par \noindent
{\bf - Step 3:} $\;$ {\it Inside $\Omega$:}
\par \noindent
We consider now $\Omega_{\varepsilon'}:=\{x\in \Omega, \, dist(x, \partial \Omega) > \varepsilon'\}.$
Set $$B:= \min_{\Omega_{\varepsilon'}} \varphi_1(x)>0.$$
We have in $\Omega_{\varepsilon'}$ by (\ref{10})  and (\ref{11})
$$z(x) \, =\, \frac{\alpha}{\lambda_1 - \mu} \varphi_1(x)  \, + \, z^{\bot}(x)  \leq 
 \frac{\alpha}{\lambda_1 - \mu} B + C_2\|h^\bot\|_{L^q} <0$$
if we choose
$$\mu - \lambda_1 < \frac{\alpha \min( B, A/2) }{C_2 \|h^\bot\|_{L^q}}. $$ 
We derive now Theorem 1. $\bullet$
\subsection{ An example}
\par \noindent
Let $N=1$, $\Omega = ]0,1[$ and $h=h_1 \varphi_1 + h_2 \varphi_2$
with $h_1>0$, $h_2>0$. We note that
\begin{equation}\label{16} \varphi_1(x) - s \varphi_2(x) = sin \pi x( 1 - 2 s cos \pi x) >0 \end{equation} 
in $\Omega$ implies $s \leq 1/2$.
For this example, taking $\mu = \lambda_1+\varepsilon, \varepsilon>0$, we have:
$$ z = \frac{h_1}{\lambda_1 - \mu}\varphi_1 + \frac{h_2}{\lambda_2 - \mu}\varphi_2
=- \frac{h_1}{\varepsilon} \left(\varphi_1 -\frac{\varepsilon h_2 }{h_1 (\lambda_2 -\lambda_1 - \varepsilon)}\varphi_2 \right).$$
If the Antimaximum Principle holds, $z<0$ in $\Omega$, and by ( \ref{16}), we have
$$\frac{\varepsilon h_2 }{h_1(\lambda_2 - \lambda_1 - \varepsilon)} \leq \frac{1}{2},$$
hence
$$\varepsilon \leq \frac{h_1(\lambda_2 - \lambda_1)}{2 h_2( 1+\frac{h_1}{2 h_2}) }
 \leq \frac{h_1(\lambda_2 - \lambda_1)}{2 h_2 }.
$$
We obtain an estimate of $\delta(h)$ similar to that  in Theorem 1.

\section{A non-cooperative system}
Now we will consider the $2 \times 2$ non-cooperative system depending on  a real parameter $\mu$:
$$-\Delta u \, =\, au\,+\, bv \, +\, \mu u \,  +\, f\;\, {\rm in }\; \Omega ,   \eqno(S_1)$$ 
$$-\Delta v \, =\,  cu \, + \, dv \, +\, \mu v \, +\, g \;\, {\rm in }\; \Omega,   \eqno(S_2)$$
$$ u=v=0 \; {\rm on } \; \partial \Omega.  \eqno(S_3)$$ 
or shortly $$-\Delta U = AU + \mu U + F \; {\rm in }  \, \Omega \, ,  \; U=0 \, {\rm on}  \; \partial \Omega . \eqno(S)$$
\par \noindent 
{\bf Hypothesis $(H_1)$} $\; $ We assume $b > 0 \,, c < 0,  \, $ and  \begin{equation} \label{D}
 D : = (a-d)^2 + 4bc >0.
\end{equation}

\par \noindent
\subsection{Eigenvalues of the system}
As usual we say that $\mu$ is an  eigenvalue of System $(S)$ if 
$(S_1) - (S_3)$ has a non trivial solution $U = (u,v) \not = 0$ for $F\equiv 0$ and we say that 
$\mu$ is a principal  eigenvalue of System $(S)$ if there exists  $U = (u,v)$ with  $u>0, v>0$  solution to $(S)$ with 
$F \equiv 0$. 
\par \noindent
Notice that, since $(S)$ is not cooperative, it is not necessarily true that there is a lowest  principal eigenvalue $\mu_1$ and that the maximum principle holds if and only if $\mu_1>0$ (Amann \cite{Am2005}).
\par \noindent
We seek solutions $u=p\varphi_1$, $v=q\varphi_1$ to the eigenvalue problem where, as above,  $(\lambda_1, \varphi_1)$ is the principal eigenpair  for $-\Delta$ with  Dirichlet boundary conditions. 
\par \noindent
Principal eigenvalues correspond to solutions with $p,q >0$. 
The associated linear system is 
$$ ( a+ \mu - \lambda_1) p\,  +\,  bq \, = \, 0,$$
$$ cp \, + \, (d+ \mu - \lambda_1) q \, = \, 0,$$
and it follows from $(H_1)$ that $(a+ \mu - \lambda_1)$ and 
$(d+ \mu - \lambda_1)$ should have opposite signs. 
We should have 
$$Det (A+ (\mu - \lambda_1) I) = ( a+ \mu - \lambda_1)(d+ \mu - \lambda_1)-bc=0,$$
which implies by $(H_1)$ that the condition on signs is satisfied and this whatever the sign of $\mu$ could be. 
(Notice that $D>0$ implies that both roots are real and that $D=0$ gives a real double root). 
\par \noindent
We have then shown directly that our system has (at least) two principal eigenvalues. Their  signs will depend on the coefficients. If, for example, 
$a<\lambda_1$, $d<\lambda_1$, the largest one is positive. We will denote the two principal eigenvalues by 
$\mu^-_1 $ and $ \mu_1^+$ where 
\begin{equation} \label{mu} \mu_1^- \,: =\,  \lambda_1 - \xi_1 \,< \,  
  \mu_1^+ \,:=\, \lambda_1 - \xi_2 , \end{equation}
where the eigenvalues of Matrix $A$ are: 
$$\xi_1=  \frac{a+d + \sqrt D }{2} \,>\, \xi_2= \frac{a+d - \sqrt D }{2} .$$
\begin{rem} Usually the Maximum Principle holds if and only if the first eigenvalue is positive. Here by replacing 
$-\Delta$ by $-\Delta + K$ with $K>0$ large enough we may get $\mu_1^- >0$.
Nevertheless the maximum principle needs an additional condition (see Theorem \ref{thm6} and its remark).
\end{rem}

\subsection{Main Theorems}
\subsubsection{The case $\mu_1^- \, <  \,\mu  \,<  \,\mu_1^+$}
 We assume in this subsection  that the parameter $\mu$ satisfies: 
\vskip0.15cm  \noindent
 $(H_2)$ $\; $  $\mu_1^- \, <  \,\mu  \,<  \,\mu_1^+$ . 
\vskip0.15cm \noindent
\begin{thm}\label{thm2}   $\;$  Assume $(H_1), \, (H_2), \, $ and 
$$  d < a,   \leqno(H_3)$$
$$ f\geq 0, \, g\geq 0, \, f, g \not \equiv 0 , f,g \in L^q, \,q>N \, if \, N\geq 2 \, ;  \; q=2 \, if \, N=1.  \leqno(H_4)$$
Then there exists $\delta>0$, independent of $\mu$, such that if 
$$ \mu < \mu_1^- + \delta,  \leqno(H_5)$$
we get  
$$ u<0, \, v>0 \; in \; \Omega; \; \frac{\partial u}{\partial \nu} >0, \, \frac{\partial v}{\partial \nu}<0 \; on \; \partial \Omega.$$ \end{thm}

\begin{rem}\label{thm3}  $ \; $ If in the theorem above we reverse  signs of $f,g,u,v $ 
that is  $ f\leq 0, \, g\leq 0, \, f, g \not \equiv 0 $, 
then for $\mu$ satisfying $(H_5)$, we get
$$ u>0, \, v<0 \; in \; \Omega; \; \frac{\partial u}{\partial \nu} <0, \, \frac{\partial v}{\partial \nu}>0 \; on \; \partial \Omega.$$ 
Note that the counterexample in subsection $(3.3)$  shows that for $f,g$ of opposite sign( $fg<0$), $u$ or $v$ may  change sign. \end{rem}

\begin{thm}\label{thm4} 
 $\;$ Assume $(H_1), \, (H_2), \, $ and 
$$   a<d,  \leqno(H'_3)$$
$$ f\leq 0, \, g\geq 0, \, f, g \not \equiv 0 \,, f,g \in L^q, \,q>N \, if \, N\geq 2 \, ;  \; q=2 \, if \, N=1.  \leqno(H'_4)$$
Then there exists $\delta>0$, independent of $\mu$, such that if 
$$ \mu < \mu_1^-  + \delta,  \leqno(H_5)$$
we obtain 
$$ u<0, \, v<0 \; in \; \Omega; \; \frac{\partial u}{\partial \nu} >0, \, \frac{\partial v}{\partial \nu}>0 \; on \; \partial \Omega.$$ 
\end{thm}

\begin{rem}\label{thm5} 
If in the theorem above we reverse  signs of $f,g,u,v $ 
that is $ f\geq 0, \, g\leq 0, \, f, g \not \equiv 0 , $
then for $\mu$ satisfying $(H_5)$, we get
$$ u>0, \, v>0 \; in \; \Omega; \; \frac{\partial u}{\partial \nu} <0, \, \frac{\partial v}{\partial \nu}<0 \; on \; \partial \Omega.$$ 
Note that, by  the changes used in the proof of the theorem  above, the counterexample in subsection $(3.3)$  shows that for $f,g$ with same sign ($fg>0$), $u$ or $v$ may  change sign.
\end{rem}

\subsubsection{The case  $\mu  \,<  \,\mu_1^-$  }

We assume in this Section  that the parameter $\mu$ satisfies: 
\vskip0.15cm  \noindent
$(H'_2)$ $\; $  $\mu  \,<  \,\mu_1^-$ . 
\vskip0.15cm \noindent
\begin{thm} \label{thm6}  $\;$  Assume $(H_1), \, (H'_2), \, $ and 
$$  a<d,   \leqno(H'_3)$$
$$ f\geq 0, \,  g \geq 0, \, f, g \not \equiv 0, \, f,g \in L^2. \leqno(H''_4)$$
Assume also  $t ^* g -f \geq 0, \,  t ^* g -f \not \equiv 0$  with $$t^*= \frac{d-a+\sqrt D}{-2c}.$$ 
Then  
$$ u>0, \, v>0 \; in \; \Omega; \; \frac{\partial u}{\partial \nu} <0, \, \frac{\partial v}{\partial \nu}<0 \; on \; \partial \Omega.$$ 
\end{thm}
\begin{rem}\label{thm7} 
As above we can reverse  signs of $f,g,u,v $ .
\end{rem}

\vskip0.2cm \noindent
\subsection{Counterexample: $a>d$}

We consider the system in 1 dimension
$$-u" \, =\, 4u\,+\, v \, +\, \mu u \,  +\, f\;\, {\rm in }\;  I:=]0; \pi[ ,   $$ 
$$-v" \, =\,  -u \, + \, v \, +\, \mu v \, +\, g \;\, {\rm in }\; I,   $$
$$ u(0)=u(\pi)=v(0) = v(\pi)=0 .  $$ 
\par \noindent
$\lambda_1=1$ and $\lambda_2=4$; $\varphi_1 = \sin x$, $\varphi_2 = \sin 2x$.
We compute
$\mu_1^-= 1- \frac{5 + \sqrt 5}{2}$.  
Choose $f = \varphi_1 - \frac{1}{2} \varphi_2 \, \geq 0$  and $g = kf$ with $k \neq 0$ to be determined later.
We obtain
$$u = u_1 \varphi_1 + u_2 \varphi_2 \;  {\rm  and } \; v = v_1 \varphi_1 + v_2 \varphi_2, $$
where 
$$u_1= \frac{k-\mu}{\mu^2 +3 \mu +1} \, , \; u_2= \frac{\mu - k - 3}{2(\mu^2 -3 \mu + 1 )}, $$
\noindent
1/ Choosing $\mu = -3 < \mu_1^-$, we get $v_1 =-1$ and 
$v_2 = \frac{1-3k}{38}$.
Therefore
$$ -v = \varphi_1 + \frac{3k-1}{38} \varphi_2, $$
and for $\frac{3k-1}{38} > \frac{1}{2}$, $v$ changes sign.
Hence Maximum Principle does not hold. 
\vskip0.3cm  \noindent
2/ Choosing $\mu^-_1< \mu = \mu_1^- + \epsilon$, $k=\mu_1^- + \epsilon^2$,  we have
$$\frac{u_2}{u_1} =\left(\frac{\mu - k -3}{k-\mu} \right) \left(\frac{\mu^2 + 3 \mu + 1}{ 2(\mu^2 -3 \mu + 1 )}, \right) = \left(\frac{3+\epsilon}{\epsilon}\right) 
\left(\frac{\sqrt 5 - \epsilon}{(9 + 3 \sqrt 5)  - ( 6 + \sqrt 5) \epsilon + \epsilon^2}\right).  $$
So that $\frac{u_2}{u_1} \rightarrow \infty$ as $\epsilon \rightarrow 0$. Hence
for these $f>0$, $g<0$, $u$ changes sign. 
$\bullet$. 
\subsection{Proofs of the main results}
\subsubsection{Some computations  and associate equation}
In the following we introduce 
\begin{equation} \label{gamma1} \gamma_1\, = \, \frac{1}{2}(a+d + 2 \mu - \sqrt D) \, = \, \lambda_1 + \mu - \mu_1^+; \end{equation}
\begin{equation} \label{gamma2}
\gamma_2=\frac{1}{2}(a+d +2\mu + \sqrt D) \, = \, \lambda_1 + \mu - \mu_1^-, 
\end{equation}
and some auxiliary results used in the proofs of our results. 
\begin{lem}\label{L.} $\;$ We have
$$ \mu < \mu_1^+ \; \Leftrightarrow \; \gamma_1\, < \lambda_1. \leqno(L1)$$
$$  \mu_1^- < \mu  \; \Leftrightarrow \;  \lambda_1 <\gamma_2 .\leqno(L2)$$
$$ \sqrt D < a-d  \; \Leftrightarrow \; d+ \mu < \gamma_1 < \gamma_2 < a+ \mu.\leqno(L3)$$
$$ \sqrt D < d-a  \; \Leftrightarrow \; a+ \mu < \gamma_1 < \gamma_2 < d+ \mu. \leqno(L4)$$
$$\mu < \mu_1^+ + \delta \; \Leftrightarrow \; \gamma_1 < \lambda_1+ \delta.\leqno(L5)$$
$$ \mu < \mu_1^- + \delta  \, \Leftrightarrow \, \gamma_2 < \lambda_1 + \delta . \leqno(L6)$$

\end{lem}

\subsubsection{Proofs of Theorems 2 and 3 } 
{\bf Proof of Theorem 2, $a>d$:}
\par \noindent
We introduce now
\begin{equation} \label{w}  w\, =\, u+tv,  \end{equation} with 
\begin{equation} \label{t}t= 
 \frac{a-d+\sqrt D}{-2c}\,  =\,  \frac{2b}{a-d-\sqrt D}
\end{equation}

so that 
\begin{equation} - \Delta w = \gamma_1 w + f+tg \, in \, \Omega; \, w|_{\partial \Omega} = 0.
 \end{equation}
We remark that
\begin{equation}
t= \frac {b}{\gamma_1 - d -\mu} = \frac{b}{a +\mu - \gamma_2}  = \frac{ \gamma_1 -a- \mu}{c}  = \frac{  d +\mu - \gamma_2 }{c}.
\end{equation}
Note first that Hypothesis $(H_3)$ implies $t>0$ and $a-d>\sqrt D$. 
By $(H_2)$, $(H_4)$, and $(L1)$ in  Lemma \ref{L.},  $ \gamma_1\, < \lambda_1$,  and 
we apply the Maximum Principle
which gives 
$w>0$ on $\Omega$ and $\frac{\partial w}{\partial \nu} <0$ on  $\partial \Omega.$
We compute
 \begin{equation} a + \mu - \frac{b}{t} = a  + d +2\mu - \gamma_1 = \gamma_2, \end{equation} and 
since $v= (w-u)/t$, we derive 
$$-\Delta u = (a + \mu - \frac{b}{t} )u + \frac{b}{t} w + f = \gamma_2 u + \frac{b}{t} w + f, $$
where  $\frac{b}{t} w + f > 0$.
From  $(H_5)$ and $(L6)$, 
$\gamma_2 \leq \lambda_1 + \delta_1  $, 
where
\begin{equation} \delta_1 :=  \delta(\frac{b}{t} w + f) , \end{equation}
 we deduce from the Antimaximum Principle that
$ u<0$ on $\Omega$ and $\frac{\partial u}{\partial \nu} >0$ on  $\partial \Omega.$
Hence $cu+g>0$.
\par \noindent
Now $(H_2)$, $(L_1)$ and $(L_3)$ imply $d+\mu < \gamma_1 <\lambda_1$ and  the Maximum Principle applied to $(S_2)$ gives
$v>0$ on $\Omega$ and $\frac{\partial v}{\partial \nu} <0$ on  $\partial \Omega$.
\par \noindent
We apply now Section 1 to estimate $\delta_1$. 
\begin{equation}\label{eq:h}  h:=\frac{b}{t} w + f = (\gamma_1 - d - \mu) w+f = \sigma  \varphi_1 + h^{\perp}.\end{equation}
First we compute $\sigma$:Here we show that this is not the case 
 for non-cooperative systems (with maybe $\mu_1^- <0$).
 \par \noindent
 In this paper we use ideas concerning the Anti-Maximum Principle due to Cl\'ement and Peletier \cite{ClPe}
  (see also \cite{FGoTaT}) in order to study non-cooperative $2 \times 2$ systems. In Section 2 
we obtain precise estimates concerning the validity interval for the Anti-maximum Principle for 
one equation. We include an example.
\par \noindent
In Section 3, we consider a non-cooperative $2 \times 2$ system with constant coefficients 
depending on a real parameter $\mu$ having two real principal eigenvalues $\mu_1^- < \mu_1^+$. We obtain 
some theorems concerning the behavior of different cones of couples of functions  having positivity 
(or negativity) properties. We give several results of this type for values of $\mu_1^- < \mu$ but 
close to $\mu_1^-$ by combining the usual Maximum Principle  and the results for the Anti-Maximum Principle 
in Section 2. We actually prove only one of such theorems, all the others are proved 
just by making suitable changes of variables.  A (partial, under an additional assumption) Maximum  Principle for $\mu < \mu_1^-$ is also obtained.

\par \noindent 
Set $f=\alpha \varphi_1 + f^{\perp}$,  $g=\beta \varphi_1 + g^{\perp}$,  $w=\kappa \varphi_1 + w^{\perp}$.
Since
$$-\Delta w = \gamma_1 w + f + \frac{b}{\gamma_1 - d - \mu} g, $$
we calculate:
$$\sigma = \alpha + (\gamma_1 -d - \mu) \kappa = \alpha  \frac{\lambda_1 - d - \mu}{\lambda_1-\gamma_1}
+ \beta  \frac{b}{\lambda_1-\gamma_1}.$$
Now we estimate $\|h^{\perp}\|_{L^2}$. 
$$-\Delta w^{\perp} = \gamma_1 w^{\perp} + f^{\perp} + \frac{b}{\gamma_1 - d - \mu} g^{\perp}.$$
The variational characterization of $\lambda_2$ gives 
$$(\lambda_2 - \gamma_1) \|w^{\perp}\|_{L^2} \leq \|f^{\perp} \|_{L^2} + \frac{b}{\gamma_1 -d -\mu } \| g^{\perp}\|_{L^2}.$$
We derive from ( \ref{eq:h})
$$\|h^{\perp}\|_{L^2} \leq \|f^{\perp} \|_{L^2} +  (\gamma_1 -d  - \mu )\| w^{\perp}\|_{L^2} 
\leq  \frac{\lambda_2 - d - \mu}{\lambda_2 - \gamma_1}\|f^{\perp} \|_{L^2} 
  +  \frac{b}{\lambda_2 - \gamma_1 } \|   g^{\perp}\|_{L^2}.$$
Reasoning as in Lemma 2.1, we show that there exists a constant $C_3$ such that 
\begin{equation} \label{2.5} \|h^{\perp}\|_{L^q} \leq C_3\left(   \frac{\lambda_2 - d - \mu}{\lambda_2-\gamma_1} \|f^{\perp} \|_{L^q} +  
\frac{b}{\lambda_2 - \gamma_1  } \|g^{\perp}\|_{L^q}\right). \end{equation}
In fact for proving  (\ref{2.5}) we use the same sequence than that in Lemma 2.1 and 
we show by induction that
$$ \|z^{\perp}\|_{L^{p_j}}\, \leq \, K(j) \left( \|f^{\perp}\|_{L^{q}}+ \|g^{\perp}\|_{L^{q}}\right).$$
Now we apply the antimaximum principle to the equation
$$ -\Delta u = \gamma_2  u + h.$$
This  is possible since by $(L6)$ in Lemma \ref{L.}, $\lambda_1 < \gamma_2 < \lambda_1 + \delta_2 = \lambda_1 + \delta(h)$ where, as in Theorem 1, 
$\delta(h) = \frac{K \sigma}{\|h^{\perp}\|_{L^{q}}}$. 
\par \noindent
Moreover we notice that $\lambda_1 - \gamma_1 = \mu_1^+ - \mu \leq  \mu_1^+ - \mu_1^-$ and therefore, since $\alpha>0$  and $\beta >0$ by $(H_4)$, 
$$\sigma =  \alpha  \frac{\lambda_1 - d - \mu }{\lambda_1-\gamma_1}
+ \beta  \frac{b}{\lambda_1-\gamma_1} \, \geq\
 \, {\cal A}:= \alpha  \frac{\lambda_1 - d - \mu_1^+ }{\mu_1^+ - \mu_1^-}
+ \beta  \frac{b}{\mu_1^+ - \mu_1^-},  $$ 
and from (\ref{2.5}), we obtain
$$ \|h^{\perp}\|_{L^q} \leq {\cal B}:=  C_3\left(   \frac{\lambda_2 - d - \mu_1^-}{\lambda_2-\lambda_1} \|f^{\perp} \|_{L^q} +  
\frac{b}{\lambda_2 - \lambda_1  } \|g^{\perp}\|_{L^q}\right). $$
From the computation above we can choose
$\delta_2= \frac {K {\cal A}}{{\cal B}}$ which does not depend on $\mu$,  and the result follows.
$\bullet$
\vskip0.2cm \noindent
{\bf Proof of Theorem \ref{thm4}:$\; a<d$. \,  } We deduce this theorem from Theorem \ref{thm2} by change of variables. 
Set $\hat a = d,\,  \hat d = a$ , $\hat u = v$,  $\hat v = -u$ and $\hat f = g$ , $\hat g =-f$.
$\hat f \geq 0$, $\hat g \geq 0$, imply $\hat u <0$, $ \hat v >0$. 
 We get Theorem \ref{thm4}.
$\bullet$
 
\subsubsection{Proof of Theorem \ref{thm6} }
Since $a<d$, we have 
$t^*= \frac{d-a + \sqrt D}{-2c}  >0$. With now the change of variable $ w\, =\,  -u+t^*v,$
 as in \cite{CaMi} (see also \cite{Le}) , we can write  the system as
\begin{equation} \label{E1}-\Delta u = \gamma_1 u + (b/t^*)w+ f \; in \, \Omega,  \end{equation}
\begin{equation} \label{E2}-\Delta v = \gamma_1 v  - cw +g \; in \, \Omega  \end{equation}
\begin{equation} \label{E3}-\Delta w = \gamma_2 w + (t^*g - f) \; in \, \Omega,  \end{equation}
$$ u=v=w=0   \; on \; \partial \Omega.$$
Now   $\mu  \,<  \,\mu_1^-$, and it follows from $(L2)$ in Lemma \ref{L.} that 
$\gamma_1 < \gamma_2 < \lambda_1$. From (\ref{E3})  it follows from the Maximum Principle that $w>0$. Then 
in (\ref{E2})  $ - cw +g>0$, and again by the Maximum Principle $v>0$. Finally, since $(b/t^*)w+ f>0$ in (\ref{E1}), again 
by the Maximum Principle $u>0$. 
$\bullet$.

\vskip0.3cm \noindent
{\bf Acknowledgements}
The authors thank the referee for useful comments.
\par \noindent
J.Hern\'andez is partially supported by the project MTM2011-26119 of the DGISPI (Spain).


\begin{thebibliography}{99}
\bibitem{Am} H.Amann, {\it Fixed point equations and nonlinear eigenvalue problems in ordered Banach spaces}, {\bf SIAM Re. 18}, 4, 1976, p.620-709.
%
\bibitem{Am2005} H.Amann, {\it Maximum Principles and Principal Eigenvalues.} \underline{Ten Mathe}
 \underline{-matical Essays on Approximation in Analysis and Topology},
J. Ferrera, J. L\'opez-G\'omez, F.R. Ruíz del Portal ed., 
Elsevier, 2005, 1 - 60.  

\bibitem{ArGa} D. Arcoya, J. Gamez,  {\it Bifurcation theory and related problems: anti-maximum principle and resonance}, {\bf Comm. Part. Diff. Equat., 26}, 2001, p.1879-1911.
%

\bibitem{CaMi} G. Caristi, E. Mitidieri, {\it Maximum principles for a class of non-cooperative elliptic systems}, {\bf Delft Progress Rep. 14 }, 1990, p.33-56.
%
\bibitem{ClPe} P. Cl\'ement, L. Peletier, {\it An anti-maximum principle for second order elliptic operators.}, {\bf J. Diff. Equ. 34 }, 1979, p.218-229.
%
\bibitem{dFMi1986} D.G.de Figueiredo, E.Mitidieri , {\it A Maximum Principle for an Elliptic System and Applications to semilinear Problems},  SIAM J. Math and Anal. N17 (1986), 836-849.
%
\bibitem{FiMi1} D.G. de Figueiredo, E. Mitidieri, {\it Maximum principles for cooperative elliptic systems}, {\bf C. R. Acad. Sci. Paris 310}, 1990,  p.49-52.
%
\bibitem{FiMi2} D.G. de Figueiredo, E. Mitidieri, {\it Maximum principles for linear elliptic systems}, {\bf Quaterno Mat.   177 }, Trieste, 1988.
%
\bibitem{FGoTaT} J. Fleckinger, J. P. Gossez, P. Tak\'{a}c, F. de Th\'elin, {\it Existence, nonexistence et principe de l'antimaximum pour le p-laplacien},  {\bf C. R. Acad. Sci. Paris 321}, 1995, p.731-734.
%
\bibitem{FHeT}   J. Fleckinger, J. Hern\'andez, F. de Th\'elin, {\it On maximum principles and existence of positive solutions for some cooperative elliptic systems}, {\bf Diff. Int. Eq. 8 }, 1, 1995, p.69-85.
%
\bibitem{Le}  M.H.L\'ecureux, {\it Au-del\` a du principe du maximum pour des syst\` emes d'op\'erateurs elliptiques. } Th\` ese, Universit\'e de Toulouse, Toulouse 1, 13 juin 2008.
\bibitem{PW} M.H.Protter, H.Weinberger,  \underline{ Maximum Principles in Differential Equa} \underline{-tions}, 
 Springer-Verlag, 1984.
\end{thebibliography}
\end{document}